\RequirePackage{fix-cm}
\documentclass[10pt,letterpaper]{article}

\usepackage[
    textwidth=30pc,
    textheight=51pc,
    left=110.77pt,
    top=73.27pt,
    headheight=0pt,
    headsep=0pt,
    footskip=28pt
]{geometry}
\usepackage{times}
\usepackage{mathtools}
\usepackage{amssymb,amsthm}
\usepackage{enumitem}
\usepackage{aliascnt}  
\usepackage{xcolor}
\usepackage[hyphens]{url}
\usepackage{titlesec}
\usepackage[hidelinks]{hyperref}
\usepackage{cleveref}  
\crefname{enumi}{condition}{conditions}   
\Crefname{enumi}{Condition}{Conditions}

\makeatletter
\renewcommand{\normalsize}{%
    \@setfontsize\normalsize{10.5}{12}%
    \abovedisplayskip=11pt plus 3pt minus 3pt
    \belowdisplayskip=\abovedisplayskip
    \abovedisplayshortskip=\abovedisplayskip
    \belowdisplayshortskip=\abovedisplayskip
}
\makeatother
\DeclareMathSizes{10.5}{10.5}{7}{6}
\renewcommand{\footnotesize}{\fontsize{8.2}{10}\selectfont}
\normalsize
\titleformat{\section}{\normalfont\normalsize\bfseries}{\thesection.}{.5em}{}
\titleformat{\subsection}{\normalfont\normalsize\bfseries}{\thesubsection.}{.5em}{}
\titleformat{\subsubsection}{\normalfont\normalsize\itshape}{\thesubsubsection.}{.5em}{}
\titlespacing*{\section}{0pt}{12pt plus 2pt minus 2pt}{6pt}
\titlespacing*{\subsection}{0pt}{12pt plus 2pt minus 2pt}{6pt}
\titlespacing*{\subsubsection}{0pt}{12pt plus 2pt minus 2pt}{6pt}

\hypersetup{
    pdftitle={The logic of floor division},
    pdfauthor={Benjamin Przybocki},
    pdfkeywords={axiomatization, floor division, equational logic}
}

\newcommand{\N}{\mathbb{N}}
\newcommand{\Z}{\mathbb{Z}}
\newcommand{\Zp}{\mathbb{Z}_{>0}}
\newcommand{\R}{\mathbb{R}}
\DeclarePairedDelimiter{\floor}{\lfloor}{\rfloor}
\DeclareMathOperator{\lcm}{lcm}
\newcommand{\Exists}[1]{\exists\,#1.\:}
\newcommand{\Forall}[1]{\forall\,#1.\:}
\allowdisplaybreaks

\theoremstyle{plain}
\newtheorem{theorem}{Theorem}

\newaliascnt{proposition}{theorem}
\newtheorem{proposition}[proposition]{Proposition}
\aliascntresetthe{proposition}

\newaliascnt{lemma}{theorem}
\newtheorem{lemma}[lemma]{Lemma}
\aliascntresetthe{lemma}

\newaliascnt{corollary}{theorem}
\newtheorem{corollary}[corollary]{Corollary}
\aliascntresetthe{corollary}

\theoremstyle{definition}
\newaliascnt{definition}{theorem}
\newtheorem{definition}[definition]{Definition}
\aliascntresetthe{definition}

\newaliascnt{remark}{theorem}

\aliascntresetthe{remark}

\newaliascnt{question}{theorem}

\aliascntresetthe{question}

\title{The logic of floor division}
\author{Benjamin Przybocki\\
    \small Carnegie Mellon University, Pittsburgh, USA\\
    \small\texttt{bprzyboc@andrew.cmu.edu}}
\date{}

\begin{document}
\maketitle

\begin{abstract}
Floor division is the operation that divides two integers and discards the remainder. This operation satisfies some nontrivial identities, but these identities have never been systematically characterized. To rectify this, we present a sound and complete axiomatization for the equational theory of the natural numbers with addition and the operations for floor division by fixed positive integers. Our axiomatization consists of four simple schemata, which are easily seen to be sound. On the other hand, the proof that these axioms are complete requires performing some surprisingly subtle equational reasoning to put every term into a normal form.
\end{abstract}

\vspace{10px}

\noindent{\textbf{Keywords:}{ axiomatization, floor division, equational logic}} 

\vspace{10px}

\section{Introduction} \label{sec-intro}

In elementary algebra, one learns some basic identities governing addition and multiplication:
\begin{align*}
    x + y &= y + x \\
    x \cdot (y + z) &= x \cdot y + x \cdot z \\
    x + 0 &= x,
\end{align*}
and so on. Once we've listed all of the obvious identities, how do we know if there are any others?

The answer to this question is well-known and easy to prove once the relevant notions have been made precise. Terms in the language of $(\N,0,1,+,\cdot)$ satisfy the axioms of commutative semirings. Conversely, these axioms suffice to prove all true identities between such terms. Moreover, the axioms allow us to put every term into a normal form (namely, a sum of monomials). It is considerations like these that led to the development of abstract algebra and, more recently, term rewriting theory~\cite{rewriting}.

If multiplication is a fundamental arithmetic operation, then so is division, and it is therefore natural to ask what kind of identities govern addition and division. We are specifically interested in arithmetic over $\N$, so we study \emph{floor division} (also called integer division or Euclidean division), which is the operation that divides two integers and then takes their floor (equivalently, it is the result of performing Euclidean division and discarding the remainder). We denote the floor division of $x$ and $y$ by
\[
    \floor*{\frac{x}{y}}.
\]
We can then ask what identities govern addition and floor division, which is to say the equational theory of $(\N,0,1,+,\floor*{\bullet / \bullet})$. Actually, this is a bit too ambitious: in \Cref{sec-floor-div-arbitrary}, we show that this equational theory is undecidable, so there is no hope of a simple characterization of the true identities. We instead set our sights on the fragment of this equational theory where we only allow floor division by fixed positive integers; that is, we consider the equational theory of $(\N,0,1,+) \cup (\floor*{\bullet / d} : d \in \Zp)$. Our primary result is a sound and complete axiomatization for this equational theory:
\begin{theorem} \label{thm-ax-floor}
    The equational theory of $(\N,0,1,+) \cup (\floor*{\bullet / d} : d \in \Zp)$ is axiomatized by the axioms of commutative monoids together with the following additional axioms:
    \begin{align*}
        \tag{A1} \floor*{\frac{n}{d}} &\approx 0 \qquad \text{for all} \ n \in [0,d-1] \\
        \tag{A2} \floor*{\frac{x}{d}} + y &\approx \floor*{\frac{x + d \cdot y}{d}} \\
        \tag{A3} x &\approx \sum_{k=0}^{d-1} \floor*{\frac{x+k}{d}} \\
        \tag{A4} \floor*{\frac{\floor*{\frac{x}{d_1}}}{d_2}} &\approx \floor*{\frac{x}{d_1 \cdot d_2}}.
    \end{align*}
\end{theorem}

It is easy to check that these axioms are sound, meaning that they are valid identities over $\N$ (\Cref{thm-sound}). Our novelty is to show that they are complete, meaning that every valid identity over $\N$ between terms over $(0,1,+) \cup (\floor*{\bullet / d} : d \in \Zp)$ can be proven from these axioms (\Cref{cor-complete}). We prove completeness by showing that the axioms suffice to put every term into a normal form (\Cref{thm-nf}) and that distinct normal forms represent distinct functions on natural numbers. This normal form may be of independent interest, since it is not obvious that a normal form ought to exist. For example, consider the following identity:
\begin{equation} \label{eq:example}
\begin{split}
    &\floor*{\frac{x + \floor*{\frac{y}{3}} + \floor*{\frac{y+2}{3}}}{3}} + \floor*{\frac{x + 2 \cdot\floor*{\frac{y+1}{3}}}{3}} \\
    ={} &\floor*{\frac{x + \floor*{\frac{y}{3}} + \floor*{\frac{y+1}{3}}}{3}} + \floor*{\frac{x + \floor*{\frac{y+1}{3}} + \floor*{\frac{y+2}{3}}}{3}},
\end{split}
\end{equation}
which is valid for all $x,y \in \N$. It is not even obvious that \eqref{eq:example} is valid, much less which side of the equation is ``simpler''.\footnote{It turns out that the right-hand side is the normal form. What makes it ``simpler'' is that it satisfies \cref{nf-cond-3} in \Cref{def-nf}.} It follows from \Cref{thm-ax-floor} that there is a derivation of \eqref{eq:example} from the axioms A1--A4, but the shortest derivation we know of has over 100 applications of A1--A4!

\subsection{Related work}

One of our axiom schemata is a variant of Hermite's identity~\cite{hermite}. Hermite's identity is usually stated as
\[
    \floor*{d \cdot x} = \sum_{k=0}^{d-1} \floor*{x + \frac{k}{d}}
\]
for all $x \in \R$ and $d \in \Zp$. Our axiom schema is
\[
    x \approx \sum_{k=0}^{d-1} \floor*{\frac{x+k}{d}}
\]
for all $d \in \Zp$. This equation is valid for all $x \in \Z$, a fact that follows from Hermite's identity. Hermite's paper is from 1884, so the study of identities involving floor division has a long history.

In 1918, Ramanujan~\cite{ramanujan} posed some problems about the floor operation, one of which was to prove
\[
    \floor*{\frac{n}{3}} + \floor*{\frac{n+2}{6}} + \floor*{\frac{n+4}{6}} = \floor*{\frac{n}{2}} + \floor*{\frac{n+3}{6}}.
\]
There are several ways to verify this, but our normal form gives a particularly perspicuous way to solve problems like this one. Namely, the equations
\begin{align*}
    \floor*{\frac{n}{3}} &= \floor*{\frac{n}{6}} + \floor*{\frac{n+3}{6}} \\
    \floor*{\frac{n}{2}} &= \floor*{\frac{n}{6}} + \floor*{\frac{n+2}{6}} + \floor*{\frac{n+4}{6}}
\end{align*}
allow us to rewrite both sides of the equation into the following normal form:
\[
    \floor*{\frac{n}{6}} + \floor*{\frac{n+2}{6}} + \floor*{\frac{n+3}{6}} + \floor*{\frac{n+4}{6}}.
\]

Bond and Levine~\cite{abelian-networks-1,abelian-networks-2,abelian-networks-3} introduced \emph{abelian processors}, which are automata whose outputs are independent of the order of their inputs. This is motivated both by models in statistical physics and by modeling distributed networks that are robust to communication delays that affect the order in which inputs are received. Holroyd, Levine, and Winkler~\cite{abelian-gates} characterized the functions $\N^m \to \N$ computable by various kinds of abelian processors. For one class of abelian processors, the computable functions are exactly the (interpretations of) terms $t$ in the language $(0,1,+) \cup (\floor*{\bullet / d} : d \in \Zp)$ such that $t(\vec{0}) = 0$. Using their characterization, they proved that every such term represents the same function as one in which the nesting depth of floor division operations is at most $m$. This also follows immediately from our normal form, which gives an independent proof of this fact.

\section{Preliminaries} \label{sec-prelim}

In this section, we discuss some notations and concepts necessary to understand the remainder of this paper.

\subsection{Notation}

Let $\Zp$ be the positive integers, $\N = \{0,1,2,\dots\}$ be the natural numbers, $\Z$ be the integers, and $\R$ be the real numbers. Given $m,n \in \Z$, let $[m,n] = \{\ell \mid m \le \ell \le n\}$ and $[n] = [1,n]$.

\subsection{Equational logic}

A \emph{signature} is a set of function symbols $\mathcal{F}$, each of which has a fixed arity. Formally, constants are nullary function symbols. \emph{Terms} are expressions built out of variables and function symbols in $\mathcal{F}$, and \emph{formulas} are expressions built out of terms together with the propositional connectives and quantifiers. A \emph{sentence} is a formula with no free variables. The only predicate symbol we use is equality, $\approx$.

An \emph{interpretation} $\mathcal{D}$ is a nonempty set $D$ together with a function $f^\mathcal{D}$ on $D$ of the appropriate arity for each $f \in \mathcal{F}$. An interpretation naturally extends to terms, and each sentence becomes true or false with respect to an interpretation. We specify an interpretation by listing its underlying set together with its signature between parentheses; for example, $(\N,0,1,+,\cdot)$.

Note that $\approx$ is a symbol within our object language, and $=$ is a symbol we use in our metalanguage. We use $\approx$ to indicate that two terms are equivalent within a given logical theory; we use $=$ to indicate that two terms are syntactically identical or that the values of the expressions are equal for all inputs with respect to some interpretation.

A \emph{theory} is a set of sentences closed under logical entailment. A set of \emph{axioms} is a set of sentences that generates a given theory when closed under logical entailment. Given an interpretation, its \emph{theory} is the set of sentences true in the interpretation; this set is obviously closed under logical entailment, justifying the terminology.

The \emph{equational fragment} of first-order logic consists of sentences of the form
\[
    \Forall{\vec{x}} s(\vec{x}) \approx t(\vec{x}),
\]
where $s$ and $t$ are terms. We leave the quantifiers implicit. An \emph{equational theory} is a set of equational sentences closed under logical entailment.

The study of the equational fragment of first-order logic is called \emph{equational logic}. Equational logic admits a simple proof system, which formalizes the intuitive notion of equational reasoning. Let $E$ be a set of equations. Then, the inference rules are as follows:
\begin{enumerate}[label=(\roman*), itemsep=0pt]
    \item If $s \approx t \in E$, then $E \vdash s \approx t$.
    \item $E \vdash t \approx t$.
    \item If $E \vdash s \approx t$, then $E \vdash t \approx s$.
    \item If $E \vdash s \approx t$ and $E \vdash t \approx u$, then $E \vdash s \approx u$.
    \item If $E \vdash s(x_1,\dots,x_n) \approx t(x_1,\dots,x_n)$, then
    \[
        E \vdash s(\sigma(x_1),\dots,\sigma(x_n)) \approx t(\sigma(x_1),\dots,\sigma(x_n)),
    \]
    where $\sigma$ is a mapping from variables to terms.
    \item If $E \vdash s_i \approx t_i$ for all $i \in [n]$, then $E \vdash f(s_1,\dots,s_n) \approx f(t_1,\dots,t_n)$, where $f$ is an $n$-ary function symbol.
\end{enumerate}
These inference rules are clearly sound, meaning that if $E \vdash s \approx t$, then $s^\mathcal{D} = t^\mathcal{D}$ for every interpretation $\mathcal{D}$ of $E$. Conversely, Birkhoff~\cite{birkhoff} showed that these inference rules are complete: if $s^\mathcal{D} = t^\mathcal{D}$ for every interpretation $\mathcal{D}$ of a set of equations $E$, then $E \vdash s \approx t$.

\subsection{Algebra}

\begin{definition}
    A \emph{commutative monoid} is an algebraic structure over the signature $(0,+)$ satisfying the following axioms:
    \begin{align*}
        0 + x &\approx x \\
        x + y &\approx y + x \\
        x + (y + z) &\approx (x + y) + z.
    \end{align*}
\end{definition}

When working with terms that contain $+$, we treat terms that are equivalent modulo the axioms of commutative monoids as if they are the same term. This allows us to write terms like $x + y + z$ without worrying about how the sum is associated or in which order we write the summands.

We also help ourselves to some notational conveniences. If $n \in \N$, we write $n$ as an abbreviation for $\underbrace{1+\dots+1}_{n \ \text{times}}$, and we write $n \cdot x$ as an abbreviation for $\underbrace{x+\dots+x}_{n \ \text{times}}$. We also use summation notation:
\[
    \sum_{i=1}^n x_i = x_1 + \dots + x_n.
\]
These constructs can only be used when $n$ is a constant.

\section{Derived identities} \label{sec-derived}

Our goal is to prove \Cref{thm-ax-floor}. We start with a few definitions. We call a term in the language of $(0,1,+) \cup (\floor*{\bullet / d} : d \in \Zp)$ a \emph{linear floor term}. We call the axioms given in \Cref{thm-ax-floor} the axioms for \emph{commutative monoids with floor division}. We say that two linear floor terms are \emph{equivalent} if they are equal modulo the axioms for commutative monoids with floor division.

The following theorem is one direction of \Cref{thm-ax-floor}:

\begin{theorem}[Soundness] \label{thm-sound}
    The axioms for commutative monoids with floor division are valid equations over $\N$.
\end{theorem}

We omit the proof of this theorem, since it is easy. Proving the converse (completeness) is the hard part of \Cref{thm-ax-floor}, which is what the rest of this paper is devoted to. In this section, we make progress toward this by deriving some identities that will be necessary for rewriting linear floor terms into our normal form.

\subsection{Basic identities} \label{sec-basic-prop}

Our first identity allows us to cancel out a common factor shared between a denominator and a coefficient in the numerator.
\begin{proposition} \label{prop-divide}
    The following equations are valid in any commutative monoid with floor division:
    \[
        \floor*{\frac{d_1 \cdot x + y}{d_1 \cdot d_2}} \approx \floor*{\frac{x + \floor*{\frac{y}{d_1}}}{d_2}}.
    \]
\end{proposition}
\begin{proof}
    We have
    \begin{align*}
        \floor*{\frac{d_1 \cdot x + y}{d_1 \cdot d_2}} &\approx \floor*{\frac{\floor*{\frac{d_1 \cdot x + y}{d_1}}}{d_2}} \qquad \text{by A4} \\
        &\approx \floor*{\frac{x + \floor*{\frac{y}{d_1}}}{d_2}} \qquad \text{by A2}. \qedhere
    \end{align*}
\end{proof}

Our second identity allows us to inflate a denominator by expanding a quotient into a sum of quotients with a larger denominator. 
\begin{proposition} \label{prop-mult}
    The following equations are valid in any commutative monoid with floor division:
    \[
        \floor*{\frac{x}{d_1}} \approx \sum_{k=0}^{d_2 - 1} \floor*{\frac{x + d_1 \cdot k}{d_1 \cdot d_2}}.
    \]
\end{proposition}
\begin{proof}
    We have
    \begin{align*}
        \floor*{\frac{x}{d_1}} &\approx \sum_{k=0}^{d_2-1} \floor*{\frac{\floor*{\frac{x}{d_1}}+k}{d_2}} \qquad \text{by A3} \\
        &\approx \sum_{k=0}^{d_2-1} \floor*{\frac{x + d_1 \cdot k}{d_1 \cdot d_2}} \qquad \text{by \Cref{prop-divide}}. \qedhere
    \end{align*}
\end{proof}

Putting the above two identities together, we can rewrite any quotient into a term in which every numerator is monic.

\begin{proposition} \label{prop-monic}
    The following equations are valid in any commutative monoid with floor division:
    \[
        \floor*{\frac{a \cdot x + b}{d}} \approx \sum_{k=0}^{a-1} \floor*{\frac{x + \floor*{\frac{b + d \cdot k}{a}}}{d}}.
    \]
\end{proposition}
\begin{proof}
    We have
    \begin{align*}
        \floor*{\frac{a \cdot x + b}{d}} &\approx \sum_{k=0}^{a-1} \floor*{\frac{a \cdot x + b + d \cdot k}{d \cdot a}} \qquad \text{by \Cref{prop-mult}} \\
        &\approx \sum_{k=0}^{a-1} \floor*{\frac{x + \floor*{\frac{b + d \cdot k}{a}}}{d}} \qquad \text{by \Cref{prop-divide}}. \qedhere
    \end{align*}
\end{proof}

Finally, we state an identity that allows us to exchange consecutive nested denominators while also reducing the number of occurrences of $x$. This will be useful in the next subsection.
\begin{proposition} \label{prop-exchange}
    The following equations are valid in any commutative monoid with floor division:
    \[
        \floor*{\frac{x + y + \floor*{\frac{x+z}{d}}}{d+1}} \approx \floor*{\frac{x + \floor*{\frac{d \cdot y + z}{d+1}}}{d}}.
    \]
\end{proposition}
\begin{proof}
    We have
    \begin{align*}
        \floor*{\frac{x + y + \floor*{\frac{x+z}{d}}}{d+1}}
        &\approx \floor*{\frac{x + \floor*{\frac{x + d \cdot y + z}{d}}}{d+1}} \qquad \text{by A2} \\
        &\approx \floor*{\frac{\floor*{\frac{(d+1) \cdot x + d \cdot y + z}{d}}}{d+1}} \qquad \text{by A2} \\
        &\approx \floor*{\frac{(d+1) \cdot x + d \cdot y + z}{d \cdot (d+1)}} \qquad \text{by A4} \\
        &\approx \floor*{\frac{x + \floor*{\frac{d \cdot y + z}{d+1}}}{d}} \qquad \text{by \Cref{prop-divide}}. \qedhere
    \end{align*}
\end{proof}

\subsection{The key identity} \label{sec-key-identity}

We start with a definition that will be essential when we define our normal form later on.

\begin{definition}
    The \emph{depth} of a variable $x$ in a linear floor term is the maximum nesting depth of floor division operations in which $x$ appears.
\end{definition}

Then, we have the following identity, which allows us to decrease the depth of a variable $x$.
\begin{proposition} \label{prop-key-identity}
    The following equations are valid in any commutative monoid with floor division:
    \[
        \floor*{\frac{\sum_{i=1}^{d_2} \floor*{\frac{x+y_i}{d_1}}}{d_2}}
        \approx
        \floor*{\frac{x + \sum_{k=0}^{d_1-1} \floor*{\frac{\sum_{i=1}^{d_2} \floor*{\frac{y_i + k}{d_1}}}{d_2}}}{d_1}}.
    \]
\end{proposition}

We call this the \emph{key identity}, because its derivation constitutes the most intricate part of our proof. We devote the rest of this subsection to proving it.

We start with the following lemma, which is proven by iterating \Cref{prop-exchange}. It allows us to eliminate $n$ copies of $\floor*{\frac{x+z}{d}}$ from a numerator.

\begin{lemma} \label{lem-elimination}
    For every $d \in \Zp$ and $n \in \N$, there is a term $V_{d,n}(y,z)$ with no variables other than $y$ and $z$ such that
    \[
        \floor*{\frac{x + y + n \cdot \floor*{\frac{x+z}{d}}}{d+n}} \approx \floor*{\frac{x + V_{d,n}(y,z)}{d}}.
    \]
    The terms $V_{d,n}$ can be defined recursively by
    \begin{align*}
        V_{d,0}(y,z) &= y, \\
        V_{d,n}(y,z) &= V_{d,n-1}\left(\floor*{\frac{(d+n-1) \cdot y + z}{d+n}}, z\right) \qquad (n > 0).
    \end{align*}
\end{lemma}
\begin{proof}
    First, we show by induction on $n$ that $V_{d,n}(z,z) \approx z$. This is immediate for $n=0$. For $n > 0$, the induction hypothesis gives
    \[
        V_{d,n}(z,z) = V_{d,n-1}\left(\floor*{\frac{(d+n) \cdot z}{d+n}}, z\right) \approx V_{d,n-1}(z,z) \approx z.
    \]

    We now prove the lemma by induction on $n$. The base case $n = 0$ is again immediate. Suppose $n > 0$, and let
    \[
        d' = d+n-1, \qquad x' = x + (n-1) \cdot \floor*{\frac{x+z}{d}}.
    \]
    The induction hypothesis, specialized to $y = z$, gives
    \[
        \floor*{\frac{x'+z}{d'}} \approx \floor*{\frac{x + V_{d,n-1}(z,z)}{d}} \approx \floor*{\frac{x+z}{d}}.
    \]
    Therefore, we have
    \begin{align*}
        \floor*{\frac{x + y + n \cdot \floor*{\frac{x+z}{d}}}{d+n}}
        &\approx \floor*{\frac{x' + y + \floor*{\frac{x'+z}{d'}}}{d'+1}} \\
        &\approx \floor*{\frac{x' + \floor*{\frac{d' \cdot y + z}{d'+1}}}{d'}} \qquad \text{by \Cref{prop-exchange}} \\
        &\approx \floor*{\frac{x + V_{d,n-1}\left(\floor*{\frac{d' \cdot y + z}{d'+1}}, z\right)}{d}} \qquad \text{by induction} \\
        &= \floor*{\frac{x + V_{d,n}(y,z)}{d}}. \qedhere
    \end{align*}
\end{proof}

\noindent
We are now ready to prove the key identity.

\begin{proof}[Proof of \Cref{prop-key-identity}]
For $d_1,d_2 \in \Zp$, write $\mathrm{I}(d_1,d_2)$ for the instance of \Cref{prop-key-identity} with parameters $d_1$ and $d_2$. We prove that $\mathrm{I}(d_1,d_2)$ is derivable by induction on $d_2$. When $d_2 = 1$, we have
\begin{align*}
    \floor*{\frac{\floor*{\frac{x+y_1}{d_1}}}{1}}
    &\approx \floor*{\frac{x+y_1}{d_1}} \qquad \text{by A3} \\
    &\approx \floor*{\frac{x + \sum_{k=0}^{d_1-1} \floor*{\frac{y_1+k}{d_1}}}{d_1}} \qquad \text{by A3} \\
    &\approx \floor*{\frac{x + \sum_{k=0}^{d_1-1} \floor*{\frac{\floor*{\frac{y_1+k}{d_1}}}{1}}}{d_1}} \qquad \text{by A3}.
\end{align*}

Now fix $d_1 \geq 1$ and $d_2 \geq 2$, and assume $\mathrm{I}(d_1',d_2')$ has been derived for all $d_1' \geq 1$ and $d_2' < d_2$. Let $n = d_2-1$ and
\[
    s = \sum_{i=2}^{d_2} \floor*{\frac{x+y_i}{d_1}}.
\]

\subsubsection*{Step 1.}

For $0 \leq k < n$, let
\[
    z_k = \sum_{\ell=0}^{d_1-1} \floor*{\frac{\sum_{i=2}^{d_2-1} \floor*{\frac{y_i+\ell}{d_1}} + \floor*{\frac{y_{d_2}+d_1 \cdot k+\ell}{d_1}}}{n}},
\]
which is $x$-free. Then, we have
\begin{align}
    \floor*{\frac{s+k}{n}} &\approx \floor*{\frac{\sum_{i=2}^{d_2-1} \floor*{\frac{x+y_i}{d_1}} + \floor*{\frac{x+y_{d_2}+d_1 \cdot k}{d_1}}}{n}} \qquad \text{by A2} \notag \\
    &\approx \floor*{\frac{x+z_k}{d_1}} \qquad \text{by $\mathrm{I}(d_1,n)$}. \label{ki:carry-rep}
\end{align}

\subsubsection*{Step 2.}

Next, we show that for all $k < n$,
\begin{equation} \label{ki:tail-z}
    \floor*{\frac{\sum_{j=k}^{n-1} \floor*{\frac{x+z_j}{d_1}}}{n-k}} \approx \floor*{\frac{x+z_k}{d_1}}.
\end{equation}

First, note that the ground term $\sum_{\ell=0}^{n-1} \floor*{\frac{\sum_{j=k}^{n-1} \floor*{\frac{j+\ell}{n}}}{n-k}}$ is equivalent to $k$. Indeed, the inner sum equals $\min(\ell,n-k)$, so its quotient by $n-k$ is $0$ for $\ell < n-k$ and $1$ for $\ell \geq n-k$. Then, we have
\begin{align*}
    \floor*{\frac{\sum_{j=k}^{n-1} \floor*{\frac{x+z_j}{d_1}}}{n-k}}
    &\approx \floor*{\frac{\sum_{j=k}^{n-1} \floor*{\frac{s+j}{n}}}{n-k}} \qquad \text{by \eqref{ki:carry-rep}} \\
    &\approx \floor*{\frac{s + \sum_{\ell=0}^{n-1} \floor*{\frac{\sum_{j=k}^{n-1} \floor*{\frac{j+\ell}{n}}}{n-k}}}{n}} \qquad \text{by $\mathrm{I}(n,n-k)$} \\
    &\approx \floor*{\frac{s+k}{n}} \\
    &\approx \floor*{\frac{x+z_k}{d_1}} \qquad \text{by \eqref{ki:carry-rep}}.
\end{align*}

\subsubsection*{Step 3.}

For each $k \le n$, define
\[
    t_k = x + d_1 \cdot \sum_{j=k}^{n-1} \floor*{\frac{x+z_j}{d_1}}, \qquad e_k = d_1 \cdot (n-k).
\]
Thus, $t_n = x$ and $t_k = t_{k+1} + d_1 \cdot \floor*{\frac{x+z_k}{d_1}}$. Let
\[
    u_0 = y_1, \qquad u_{k+1} = V_{e_k,d_1}(u_k, z_k) \quad (0 \leq k < n).
\]
We show that for all $k < n$,
\begin{equation} \label{ki:remove-carry}
    \floor*{\frac{u_k + t_k}{e_k + d_1}} \approx \floor*{\frac{u_{k+1} + t_{k+1}}{e_k}}.
\end{equation}
We have
\begin{align}
    \floor*{\frac{t_{k+1} + z_k}{e_k}}
    &\approx \floor*{\frac{\sum_{j=k+1}^{n-1} \floor*{\frac{x+z_j}{d_1}} + \floor*{\frac{x+z_k}{d_1}}}{n-k}} \qquad \text{by \Cref{prop-divide}} \notag \\
    &= \floor*{\frac{\sum_{j=k}^{n-1} \floor*{\frac{x+z_j}{d_1}}}{n-k}} \notag \\
    &\approx \floor*{\frac{x+z_k}{d_1}} \qquad \text{by \eqref{ki:tail-z}}. \label{ki:suffix-quotient}
\end{align}
Then, we have
\begin{align*}
    \floor*{\frac{u_k + t_k}{e_k + d_1}}
    &= \floor*{\frac{t_{k+1} + u_k + d_1 \cdot \floor*{\frac{x+z_k}{d_1}}}{e_k + d_1}} \\
    &\approx \floor*{\frac{t_{k+1} + u_k + d_1 \cdot \floor*{\frac{t_{k+1}+z_k}{e_k}}}{e_k + d_1}} \qquad \text{by \eqref{ki:suffix-quotient}} \\
    &\approx \floor*{\frac{t_{k+1} + V_{e_k,d_1}(u_k,z_k)}{e_k}} \qquad \text{by \Cref{lem-elimination}} \\
    &= \floor*{\frac{u_{k+1} + t_{k+1}}{e_k}}.
\end{align*}

\subsubsection*{Step 4.}

Since $d_2 = n+1$, we have
\begin{align*}
    \floor*{\frac{\sum_{i=1}^{d_2} \floor*{\frac{x+y_i}{d_1}}}{d_2}}
    &= \floor*{\frac{\floor*{\frac{x+y_1}{d_1}} + s}{n+1}} \\
    &\approx \floor*{\frac{\floor*{\frac{x+y_1}{d_1}} + \sum_{k=0}^{n-1} \floor*{\frac{s+k}{n}}}{n+1}} \qquad \text{by A3} \\
    &\approx \floor*{\frac{\floor*{\frac{x+y_1}{d_1}} + \sum_{k=0}^{n-1} \floor*{\frac{x+z_k}{d_1}}}{n+1}} \qquad \text{by \eqref{ki:carry-rep}} \\
    &\approx \floor*{\frac{x + y_1 + d_1 \cdot \sum_{k=0}^{n-1} \floor*{\frac{x+z_k}{d_1}}}{d_1 \cdot (n+1)}} \qquad \text{by \Cref{prop-divide}} \\
    &= \floor*{\frac{u_0 + t_0}{d_1 \cdot (n+1)}}.
\end{align*}
Then, by $n$ applications of \eqref{ki:remove-carry}, we have
\[
    \floor*{\frac{u_0 + t_0}{d_1 \cdot (n+1)}} \approx \floor*{\frac{u_n + t_n}{d_1}} = \floor*{\frac{x + u_n}{d_1}}.
\]
Thus,
\begin{equation} \label{ki:lhs-un}
    \floor*{\frac{\sum_{i=1}^{d_2} \floor*{\frac{x+y_i}{d_1}}}{d_2}} \approx \floor*{\frac{x+u_n}{d_1}}.
\end{equation}
It remains to identify $u_n$. We have
\begin{align*}
    u_n &\approx \sum_{k=0}^{d_1-1} \floor*{\frac{u_n+k}{d_1}} \qquad \text{by A3} \\
    &\approx \sum_{k=0}^{d_1-1} \floor*{\frac{\sum_{i=1}^{d_2} \floor*{\frac{y_i+k}{d_1}}}{d_2}} \qquad \text{by \eqref{ki:lhs-un}}.
\end{align*}
Making this substitution for $u_n$ in \eqref{ki:lhs-un} yields $\mathrm{I}(d_1,d_2)$.
\end{proof}

\subsection{Consequences of the key identity} \label{sec-more-prop}

Now we record two consequences of the key identity. The first identity takes a quotient in which $x$ is in a depth-1 position and $y$ is in a depth-2 position and rewrites it into one in which $y$ is in a depth-1 position and $x$ is in a depth-2 position.

\begin{proposition} \label{prop-transpose}
    The following equations are valid in any commutative monoid with floor division:
    \[
        \floor*{\frac{x + \sum_{i=1}^{d_2} \floor*{\frac{y+z_i}{d_1}}}{d_2}} \approx \floor*{\frac{y + \sum_{k=0}^{d_1-1} \floor*{\frac{x + \sum_{i=1}^{d_2} \floor*{\frac{z_i + k}{d_1}}}{d_2}}}{d_1}}.
    \]
\end{proposition}
\begin{proof}
    We have
    \begin{align*}
        &\floor*{\frac{x + \sum_{i=1}^{d_2} \floor*{\frac{y+z_i}{d_1}}}{d_2}} \\
        \approx{} &\floor*{\frac{\floor*{\frac{y+d_1 \cdot x+z_1}{d_1}} + \sum_{i=2}^{d_2} \floor*{\frac{y+z_i}{d_1}}}{d_2}} \qquad \text{by A2} \\
        \approx{} &\floor*{\frac{y + \sum_{k=0}^{d_1-1} \floor*{\frac{\floor*{\frac{d_1 \cdot x + z_1 + k}{d_1}} + \sum_{i=2}^{d_2} \floor*{\frac{z_i + k}{d_1}}}{d_2}}}{d_1}} \qquad \text{by \Cref{prop-key-identity}} \\
        \approx{} &\floor*{\frac{y + \sum_{k=0}^{d_1-1} \floor*{\frac{x + \sum_{i=1}^{d_2} \floor*{\frac{z_i + k}{d_1}}}{d_2}}}{d_1}} \qquad \text{by A2}. \qedhere
    \end{align*}
\end{proof}

The next identity is one step of our normal form procedure. It can be seen as ``normalizing'' the terms $y_i$ in the left-hand side.

\begin{proposition} \label{prop-normalize-yi}
    The following equations are valid in any commutative monoid with floor division:
    \[
        \sum_{i=1}^n \floor*{\frac{x+y_i}{d}} \approx \sum_{k=0}^{n-1} \floor*{\frac{x + \sum_{\ell=0}^{d-1} \floor*{\frac{\sum_{i=1}^n \floor*{\frac{y_i+\ell}{d}} + k}{n}}}{d}}.
    \]
\end{proposition}
\begin{proof}
    We have
    \begin{align*}
        \sum_{i=1}^n \floor*{\frac{x+y_i}{d}} &\approx \sum_{k=0}^{n-1} \floor*{\frac{\sum_{i=1}^n \floor*{\frac{x+y_i}{d}} + k}{n}} \qquad \text{by A3} \\
        &\approx \sum_{k=0}^{n-1} \floor*{\frac{x + \sum_{\ell=0}^{d-1} \floor*{\frac{\sum_{i=1}^n \floor*{\frac{y_i+\ell}{d}} + k}{n}}}{d}} \qquad \text{by \Cref{prop-transpose}}. \qedhere
    \end{align*}
\end{proof}

\section{Normal form and completeness} \label{sec-norm}

In this section, we complete the proof of \Cref{thm-ax-floor}. We first define our normal form and show that distinct normal forms represent distinct functions on natural numbers. Then, we show how to rewrite every linear floor term into a normal form using the axioms for commutative monoids with floor division.

\subsection{Normal form} \label{sec-nf-def}

We write $f(x_1,\dots,x_n)$ to denote a linear floor term in the variables $x_1,\dots,x_n$ such that each $x_i$ actually appears in the term (i.e., we omit variables that do not appear in the term).
\begin{definition} \label{def-nf}
    A linear floor term is a \emph{normal form} if it is of the form
    \[
        f() = k
    \]
    for some $k \in \N$ or
    \[
        f(x_1,\dots,x_m) = k + \sum_{i=0}^{n-1} \floor*{\frac{x_1 + f_i(x_2,\dots,x_m)}{d}}
    \]
    for some $k \in \N$, $n, d \in \Zp$, and $f_i : \N^{m-1} \to \N$, where
    \begin{enumerate}[label=(\roman*), itemsep=0pt]
        \item $f_i(x_2,\dots,x_m)$ is a normal form for all $i < n$; \label{nf-cond-1}
        \item $f_i(0,\dots,0) < d$ for all $i < n$; \label{nf-cond-2}
        \item we have
        \[
            f_i(x_2,\dots,x_m) \approx \sum_{\ell=0}^{d-1} \floor*{\frac{\sum_{j=0}^{n-1} \floor*{\frac{f_j(x_2,\dots,x_m)+\ell}{d}} + i}{n}}
        \]
        for all $i < n$; and \label{nf-cond-3}
        \item there is no $d' > 1$ such that $d' \mid n$, $d' \mid d$, and
        \[
            f_{i+n/d'}(x_2,\dots,x_m) \approx f_i(x_2,\dots,x_m) + d/d'
        \]
        for all $i < n - n/d'$. \label{nf-cond-4}
    \end{enumerate}
\end{definition}

We briefly describe the purpose of each condition in \Cref{def-nf}. \Cref{nf-cond-1} makes the definition recursive: the terms $f_i$ are themselves normal forms in the remaining variables. \Cref{nf-cond-2} says that no numerator can be reduced by moving a multiple of $d$ out of its quotient and into the constant $k$ using A2; in particular, $k = f(0,\dots,0)$. \Cref{nf-cond-3} says that the terms $f_i(x_2,\dots,x_m)$ are normalized with respect to the \Cref{prop-normalize-yi} rewrite. Finally, \Cref{nf-cond-4} says that the term cannot be rewritten, using the right-to-left direction of \Cref{prop-mult}, as a sum of only $n/d'$ quotients with the smaller denominator $d/d'$; in other words, $d$ is as small as possible.

\begin{lemma} \label{lem-distinct}
    Distinct normal forms represent distinct functions on $\N^m$.
\end{lemma}
\begin{proof}
    We induct on $m$. This is clear if $m=0$, so suppose $m>0$. First, let
    \begin{align*}
        f(x_1,\dots,x_m) &= k + h(x_1,\dots,x_m), \\
        h(x_1,\dots,x_m) &= \sum_{i=0}^{n-1} \floor*{\frac{x_1 + f_i(x_2,\dots,x_m)}{d}}
    \end{align*}
    be a normal form, and let
    \[
        \Delta h(x_1,\dots,x_m) = h(x_1+1,x_2,\dots,x_m) - h(x_1,\dots,x_m).
    \]
    Call $p \in \Zp$ a \emph{period} if
    \[
        \Delta h(x_1+p,x_2,\dots,x_m) = \Delta h(x_1,\dots,x_m) \qquad \text{for all } (x_1,\dots,x_m) \in \N^m.
    \]
    We claim that $d$ is the least period. Since $h(x_1+d,x_2,\dots,x_m) = h(x_1,\dots,x_m) + n$, the number $d$ is a period. Let $p$ be the least period. Writing $d = cp + r$ with $0 \le r < p$, we have
    \[
        \Delta h(x_1 + r, x_2,\dots,x_m) = \Delta h(x_1 + r + cp, x_2,\dots,x_m) = \Delta h(x_1,\dots,x_m),
    \]
    so $r = 0$ by the minimality of $p$. Since $p$ is a period, the difference
    \[
        h(x_1+p,x_2,\dots,x_m) - h(x_1,\dots,x_m)
    \]
    does not depend on $x_1$, and summing over $c$ consecutive blocks of length $p$ shows that $c$ divides $n$ and
    \[
        h(x_1+p,x_2,\dots,x_m) = h(x_1,\dots,x_m) + n/c.
    \]
    By \cref{nf-cond-3} of \Cref{def-nf} (and soundness), for all $i < n - n/c$ we have
    \begin{align*}
        f_{i+n/c}(x_2,\dots,x_m) &= \sum_{\ell=0}^{d-1} \floor*{\frac{h(\ell,x_2,\dots,x_m) + i + n/c}{n}} \\
        &= \sum_{\ell=p}^{d+p-1} \floor*{\frac{h(\ell,x_2,\dots,x_m) + i}{n}} \\
        &= \sum_{\ell=0}^{d-1} \floor*{\frac{h(\ell,x_2,\dots,x_m) + i}{n}} + p = f_i(x_2,\dots,x_m) + p,
    \end{align*}
    where the third equality uses $h(\ell+d,x_2,\dots,x_m) = h(\ell,x_2,\dots,x_m) + n$. If $c > 1$, then $f_{i+n/c}(x_2,\dots,x_m)$ and $f_i(x_2,\dots,x_m) + p$ are normal forms (after combining constants) representing the same function, so by the induction hypothesis they are identical, and in particular equivalent. This contradicts \cref{nf-cond-4} of \Cref{def-nf} with $d' = c$. Hence $c = 1$ and $p = d$.

    Now let
    \begin{align*}
        f(x_1,\dots,x_m) &= k_1 + \sum_{i=0}^{n_1-1} \floor*{\frac{x_1 + f_i(x_2,\dots,x_m)}{d_1}} \\
        g(x_1,\dots,x_m) &= k_2 + \sum_{i=0}^{n_2-1} \floor*{\frac{x_1 + g_i(x_2,\dots,x_m)}{d_2}}
    \end{align*}
    be linear floor terms in normal form such that $f(x_1,\dots,x_m) = g(x_1,\dots,x_m)$ for all $(x_1,\dots,x_m) \in \N^m$. We have $k_1 = f(0,\dots,0) = g(0,\dots,0) = k_2$; write $k$ for this value. By the claim, $d_1$ and $d_2$ are both the least period of the same function, so $d_1 = d_2$; write $d$ for this value. Then $n_1 = f(d,0,\dots,0) - k = g(d,0,\dots,0) - k = n_2$; write $n$ for this value. By \cref{nf-cond-3} of \Cref{def-nf} (and soundness), for all $i < n$,
    \begin{align*}
        f_i(x_2,\dots,x_m) &= \sum_{\ell=0}^{d-1} \floor*{\frac{f(\ell,x_2,\dots,x_m) - k + i}{n}} \\
        &= \sum_{\ell=0}^{d-1} \floor*{\frac{g(\ell,x_2,\dots,x_m) - k + i}{n}} = g_i(x_2,\dots,x_m)
    \end{align*}
    for all $(x_2,\dots,x_m) \in \N^{m-1}$. By the induction hypothesis, $f_i$ and $g_i$ are identical as linear floor terms for all $i < n$, so $f$ and $g$ are identical.
\end{proof}

\subsection{Rewriting into normal form} \label{sec-nf-exist}

Now we show how to rewrite every linear floor term into a normal form. We start with a couple of lemmas. Our first lemma shows how to rewrite all quotients involving $x_1$ to have a common denominator.

\begin{lemma} \label{lem-x1-d1-almost-normal}
    Every linear floor term in which $x_1$ has depth 1 is equivalent to a linear floor term of the form
    \[
        k + \sum_{i=0}^{n-1} \floor*{\frac{x_1 + f_i(x_2,\dots,x_m)}{d}}
    \]
    for some $k \in \N$, $n, d \in \Zp$, and $f_i : \N^{m-1} \to \N$.
\end{lemma}
\begin{proof}
    Let $f(x_1,\dots,x_m)$ be a linear floor term in which $x_1$ has depth 1. Using \Cref{prop-monic}, we can assume that every floor division in which $x_1$ appears is of the form
    \[
        \floor*{\frac{x_1 + f_i(x_2,\dots,x_m)}{d_i}}
    \]
    for some $d_i \in \Zp$ and $f_i : \N^{m-1} \to \N$. Using A2, we can absorb any $x_1$-free additive terms into one of these quotients. Thus, we can write
    \[
        f(x_1,\dots,x_m) \approx a \cdot x_1 + k + \sum_{i=0}^{n-1} \floor*{\frac{x_1 + f_i(x_2,\dots,x_m)}{d_i}}
    \]
    for some $a,k,n \in \N$, $d_i \in \Zp$, and $f_i(x_2,\dots,x_m) : \N^{m-1} \to \N$. Using A3, we can rewrite this as
    \[
        f(x_1,\dots,x_m) \approx k + a \cdot \floor*{\frac{x_1}{1}} + \sum_{i=0}^{n-1} \floor*{\frac{x_1 + f_i(x_2,\dots,x_m)}{d_i}}.
    \]
    Let $d = \lcm(d_0, \dots, d_{n-1})$. Using \Cref{prop-mult}, we can rewrite $f(x_1,\dots,x_m)$ as
    \[
        f(x_1,\dots,x_m) \approx k + \sum_{i=0}^{n'-1} \floor*{\frac{x_1 + f'_i(x_2,\dots,x_m)}{d}}
    \]
    for some $n' \in \Zp$ and $f'_i(x_2,\dots,x_m) : \N^{m-1} \to \N$, and the right-hand side is of the desired form.
\end{proof}

\noindent
Our next lemma shows how to rewrite any linear floor term so that $x_1$ has depth 1.

\begin{lemma} \label{lem-x1-d1}
    Every linear floor term in which $x_1$ appears is equivalent to a linear floor term in which $x_1$ has depth 1.
\end{lemma}
\begin{proof}
    It suffices to show that every floor division of the form
    \[
        \floor*{\frac{f(x_1,\dots,x_m)}{d}},
    \]
    where $f(x_1,\dots,x_m)$ is a linear floor term, is equivalent to a linear floor term in which $x_1$ has depth 1. By induction, we may assume $f(x_1,\dots,x_m)$ is equivalent to a linear floor term in which $x_1$ has depth 1. By \Cref{lem-x1-d1-almost-normal}, we can write
    \[
        f(x_1,\dots,x_m) \approx k + \sum_{i=0}^{n-1} \floor*{\frac{x_1 + f_i(x_2,\dots,x_m)}{d'}}
    \]
    for some $k \in \N$, $n, d' \in \Zp$, and $f_i : \N^{m-1} \to \N$. By \Cref{prop-mult}, we can write
    \begin{align*}
        f(x_1,\dots,x_m) &\approx k + \sum_{i=0}^{n-1} \sum_{k'=0}^{d-1} \floor*{\frac{x_1 + f_i(x_2,\dots,x_m) + d' \cdot k'}{d \cdot d'}} \\
        &\approx k + \sum_{i=0}^{n \cdot d - 1} \floor*{\frac{x_1 + f'_i(x_2,\dots,x_m)}{d \cdot d'}}
    \end{align*}
    for some $f'_i : \N^{m-1} \to \N$. Similarly, we can write
    \[
        \floor*{\frac{f(x_1,\dots,x_m)}{d}} \approx \sum_{k'=0}^{n-1} \floor*{\frac{k + d \cdot k' + \sum_{i=0}^{n \cdot d - 1} \floor*{\frac{x_1 + f'_i(x_2,\dots,x_m)}{d \cdot d'}}}{n \cdot d}}.
    \]
    Then, by \Cref{prop-transpose}, we can write
    \[
        \floor*{\frac{f(x_1,\dots,x_m)}{d}} \approx \sum_{k'=0}^{n-1} \floor*{\frac{x_1 + f''_{k'}(x_2,\dots,x_m)}{d \cdot d'}}
    \]
    for some $f''_{k'}(x_2,\dots,x_m) : \N^{m-1} \to \N$, and the right-hand side is a linear floor term in which $x_1$ has depth 1.
\end{proof}

\begin{theorem} \label{thm-nf}
    Every linear floor term is equivalent to a normal form.
\end{theorem}
\begin{proof}
    This is clear if $f()$ is a linear floor term in which no variables appear, so suppose $f(x_1,\dots,x_m)$ is a linear floor term, where $m>0$. By \Cref{lem-x1-d1-almost-normal,lem-x1-d1}, $f(x_1,\dots,x_m)$ is equivalent to a linear floor term of the form
    \[
        f(x_1,\dots,x_m) \approx k + \sum_{i=0}^{n-1} \floor*{\frac{x_1 + f_i(x_2,\dots,x_m)}{d}}
    \]
    for some $k \in \N$, $n, d \in \Zp$, and $f_i : \N^{m-1} \to \N$. We proceed by induction on $m$, so we may assume that each $f_i(x_2,\dots,x_m)$ is in normal form. Using A2, we can rewrite this expression to reduce each $f_i(0,\dots,0) \ge d$, yielding a term of the form
    \[
        f(x_1,\dots,x_m) \approx k' + \sum_{i=0}^{n-1} \floor*{\frac{x_1 + f'_i(x_2,\dots,x_m)}{d}}
    \]
    for some $k' \in \N$ and $f'_i : \N^{m-1} \to \N$ satisfying $f'_i(0,\dots,0) < d$ for all $i < n$. Write
    \[
        h(x_1,\dots,x_m) = \sum_{j=0}^{n-1} \floor*{\frac{x_1 + f'_j(x_2,\dots,x_m)}{d}},
    \]
    so that $f(x_1,\dots,x_m) \approx k' + h(x_1,\dots,x_m)$ and $h(0,\dots,0) = 0$. For $i < n$, let
    \[
        f''_i(x_2,\dots,x_m) = \sum_{\ell=0}^{d-1} \floor*{\frac{h(\ell,x_2,\dots,x_m) + i}{n}}.
    \]
    By \Cref{prop-normalize-yi},
    \begin{align}
        h(x_1,\dots,x_m) &\approx \sum_{i=0}^{n-1} \floor*{\frac{x_1 + f''_i(x_2,\dots,x_m)}{d}}, \label{nf:rep} \\
        f''_i(x_2,\dots,x_m) &\approx \sum_{\ell=0}^{d-1} \floor*{\frac{h(\ell,x_2,\dots,x_m) + i}{n}} \qquad (i < n). \label{nf:recover}
    \end{align}
    Among all $n, d \in \Zp$ and $x_1$-free terms $f''_0,\dots,f''_{n-1}$ satisfying \eqref{nf:rep} and \eqref{nf:recover}, choose one with $d$ minimal.

    By the induction hypothesis, we can put each $f''_i(x_2,\dots,x_m)$ in normal form to obtain
    \[
        f(x_1,\dots,x_m) \approx k' + \sum_{i=0}^{n-1} \floor*{\frac{x_1 + f'''_i(x_2,\dots,x_m)}{d}},
    \]
    where each $f'''_i(x_2,\dots,x_m)$ is a normal form equivalent to $f''_i(x_2,\dots,x_m)$. In particular,
    \[
        h(x_1,\dots,x_m) \approx \sum_{j=0}^{n-1} \floor*{\frac{x_1 + f'''_j(x_2,\dots,x_m)}{d}}.
    \]
    We claim that this expression is in normal form. We verify the conditions from \Cref{def-nf} as follows:
    \begin{enumerate}[label=(\roman*), itemsep=0pt]
        \item By construction, $f'''_i(x_2,\dots,x_m)$ is a normal form for all $i < n$.
        \item By soundness,
        \[
            0 = h(0,\dots,0) = \sum_{i=0}^{n-1} \floor*{\frac{f'''_i(0,\dots,0)}{d}},
        \]
        so $f'''_i(0,\dots,0) < d$ for all $i < n$.
        \item For all $i < n$, we have
        \begin{align*}
            f'''_i(x_2,\dots,x_m) &\approx f''_i(x_2,\dots,x_m) \\
            &\approx \sum_{\ell=0}^{d-1} \floor*{\frac{h(\ell,x_2,\dots,x_m) + i}{n}} \qquad \text{by \eqref{nf:recover}} \\
            &\approx \sum_{\ell=0}^{d-1} \floor*{\frac{\sum_{j=0}^{n-1} \floor*{\frac{f'''_j(x_2,\dots,x_m)+\ell}{d}} + i}{n}}.
        \end{align*}
        \item Suppose for the sake of contradiction that there is a $d' > 1$ such that $d' \mid n$, $d' \mid d$, and $f'''_{i+n/d'}(x_2,\dots,x_m) \approx f'''_i(x_2,\dots,x_m) + d/d'$ for all $i < n - n/d'$. Let $N = n/d'$ and $p = d/d'$. Since $f'''_i \approx f''_i$ for all $i < n$, \eqref{nf:rep} and \eqref{nf:recover} also hold with $f'''_i$ in place of $f''_i$. Then, by our assumption and \Cref{prop-mult},
        \begin{align*}
            h(x_1,\dots,x_m) &\approx \sum_{i=0}^{N-1} \sum_{s=0}^{d'-1} \floor*{\frac{x_1 + f'''_i(x_2,\dots,x_m) + sp}{d' \cdot p}} \\
            &\approx \sum_{i=0}^{N-1} \floor*{\frac{x_1 + f'''_i(x_2,\dots,x_m)}{p}},
        \end{align*}
        so, by A2, $h(\ell + sp, x_2,\dots,x_m) \approx h(\ell,x_2,\dots,x_m) + sN$ for all $\ell, s \in \N$. Hence, for all $i < N$,
        \begin{align*}
            f'''_i(x_2,\dots,x_m) &\approx \sum_{\ell=0}^{d-1} \floor*{\frac{h(\ell,x_2,\dots,x_m) + i}{n}} \\
            &= \sum_{\ell=0}^{p-1} \sum_{s=0}^{d'-1} \floor*{\frac{h(\ell+sp,x_2,\dots,x_m) + i}{d' \cdot N}} \\
            &\approx \sum_{\ell=0}^{p-1} \sum_{s=0}^{d'-1} \floor*{\frac{h(\ell,x_2,\dots,x_m) + sN + i}{d' \cdot N}} \\
            &\approx \sum_{\ell=0}^{p-1} \floor*{\frac{h(\ell,x_2,\dots,x_m) + i}{N}} \qquad \text{by \Cref{prop-mult}}.
        \end{align*}
        Thus, \eqref{nf:rep} and \eqref{nf:recover} hold with $N$, $p$, and $f'''_i$ in place of $n$, $d$, and $f''_i$. Since $p < d$, this contradicts the minimality of $d$. \qedhere
    \end{enumerate}
\end{proof}

\subsection{Completeness} \label{sec-nf-complete}

With our normal form result in hand, it is now easy to prove that the axioms for commutative monoids with floor division are complete.

\begin{corollary}[Completeness] \label{cor-complete}
    If two linear floor terms represent the same function on $\N^m$, then they are equivalent.
\end{corollary}
\begin{proof}
    Let $f$ and $g$ be linear floor terms that represent the same function on $\N^m$. By \Cref{thm-nf}, $f$ and $g$ are each equivalent to some normal form. By \Cref{thm-sound}, the normal forms of $f$ and $g$ also represent the same function on $\N^m$, so these normal forms are identical by \Cref{lem-distinct}. Hence, $f$ and $g$ are equivalent.
\end{proof}

\begin{proof}[Proof of \Cref{thm-ax-floor}]
    This is immediate from \Cref{thm-sound} and \Cref{cor-complete}.
\end{proof}

\section*{Acknowledgments}
    I am grateful to Siddharth Bhat and Tobias Grosser for advising my master's thesis, which this paper is based on. This research was supported by the DARPA expMath program through the DARPA CMO contract number HR0011262E028, by the National Science Foundation (NSF) under grant DMS2434625, and by the NSF Graduate Research Fellowship Program under Grant Nos. DGE-2140739 and DGE-2631988.

\bibliographystyle{unsrt}
\begingroup
\footnotesize
\bibliography{bib}
\endgroup

\appendix
\crefalias{section}{appendix}

\section{Floor division by arbitrary terms} \label{sec-floor-div-arbitrary}

So far we have only considered floor division by a constant $d \in \Zp$. What if we expand our language to allow floor divisions of the form
\[
    \floor*{\frac{x}{y}},
\]
where $x$ and $y$ are arbitrary terms? It is not clear what the semantics of this expression should be when $y = 0$. We will sidestep this problem by only allowing floor divisions of the form
\[
    \floor*{\frac{x}{1+y}},
\]
which ensures that the denominator is always positive. Formally, we consider the equational theory of $\mathcal{D} \coloneqq (\N,0,1,+,\floor*{\bullet / (1 + \bullet)})$, where $\floor*{\bullet / (1 + \bullet)}$ is a binary function on $\N$ interpreted in the obvious way.

We now show that the equational theory of $\mathcal{D}$ is undecidable. Our proof uses Matiyasevich's famous negative resolution of Hilbert's tenth problem~\cite{h10}. Given a first-order theory, let its \emph{Diophantine fragment} consist of sentences of the form
\[
    \Exists{\vec{x}} \bigwedge_{i=1}^n s_i(\vec{x}) \approx t_i(\vec{x}),
\]
where $s_i$ and $t_i$ are terms. A \emph{Diophantine theory} is a set of Diophantine sentences closed under logical entailment. Matiyasevich~\cite{h10} showed that the Diophantine theory of $(\N,0,1,+,\cdot)$ is undecidable. Our proof proceeds by reducing the decision problem for the Diophantine theory of $(\N,0,1,+,\cdot)$ to the decision problem for the equational theory of $\mathcal{D}$.

\begin{lemma} \label{lem-ex-red}
    The decision problem for the Diophantine theory of $\mathcal{D}$ is Turing reducible to the decision problem for the equational theory of $\mathcal{D}$.
\end{lemma}
\begin{proof}
    Let $f_i$ and $g_i$ be terms in the language of $(0,1,+,\floor*{\bullet / (1 + \bullet)})$, for all $i \in [n]$. Then,
    \[
        \Exists{\vec{x}} \bigwedge_{i=1}^n f_i(\vec{x}) \approx g_i(\vec{x})
    \]
    is true if and only if
    \[
        \Forall{\vec{x}} \floor*{\frac{1}{1 + \sum_{i=1}^n \left(\floor*{\frac{f_i(\vec{x})}{1+g_i(\vec{x})}} + \floor*{\frac{g_i(\vec{x})}{1+f_i(\vec{x})}}\right)}} \approx 0
    \]
    is false.
\end{proof}

\begin{lemma} \label{lem-mult-def}
    Multiplication is definable in the Diophantine theory of $\mathcal{D}$.
\end{lemma}
\begin{proof}
    For all $n, x, y \in \N$, we have $n = x \cdot y$ if and only if
    \[
        y = \floor*{\frac{n+y}{1+x}} \quad \text{and} \quad y = \floor*{\frac{n+y+x}{1+x}}. \qedhere
    \]
\end{proof}

\begin{theorem} \label{thm-floor-div-undecidable}
    The equational theory of $\mathcal{D}$ is undecidable.
\end{theorem}
\begin{proof}
    By \Cref{lem-ex-red,lem-mult-def}, the decision problem for the Diophantine theory of $(\N,0,1,+,\cdot)$ is Turing reducible to the decision problem for the equational theory of $\mathcal{D}$. By Matiyasevich's theorem~\cite{h10}, the Diophantine theory of $(\N,0,1,+,\cdot)$ is undecidable.
\end{proof}

\end{document}